\documentclass[12pt]{article}
\renewcommand{\baselinestretch}{1.5}
\usepackage{epsfig}\usepackage{epstopdf}
\usepackage{pstricks}
\usepackage{amsmath}
\usepackage{latexsym}
\usepackage{amssymb}
\usepackage{natbib}
\usepackage{bm}

\newtheorem{ex}{Example}[section]
\newtheorem{thm}[ex]{Theorem}

\newenvironment{proof}{\noindent{\bf Proof:}}{$\hfill \Box$ \vspace{12pt}}

\begin{document}

\renewcommand{\baselinestretch}{1.5}

\begin{center}
{\large \textbf{Discussion: The NoMax Strategy and Correlated Outputs}}
\end{center}
\begin{center}

\begin{tabular}{cc}
Warren Hare & Jason L.\ Loeppky  \\[-7pt]
Mathematics & Statistics
\\[-7pt]
University of British Columbia & University of British Columbia\\[-7pt]
Kelowna, BC  V1V 1V7, CANADA & Kelowna, BC  V1V 1V7, CANADA
\\[-7pt]
(warren.hare@ubc.ca) & (jason.loeppky@ubc.ca)\end{tabular}
\begin{tabular}{c}
Brian J. Williams\\[-7pt]
Statistical Sciences Group
\\[-7pt]
 Los Alamos National Laboratory \\[-7pt]
Los Alamos, NM, 87545, USA \\[-7pt]

(brianw@lanl.gov)

\end{tabular}
\end{center}
\renewcommand{\baselinestretch}{1.0}

\noindent {\bf \today}

\renewcommand{\baselinestretch}{1.6}

\section{Introduction}
\label{sect:intro}

In \cite{GraGraLeD2015}, the authors are to be congratulated for providing a statistically-motivated practical perspective to the important problem of optimizing a black-box function over a constrained region.  They assume that the constraint sets are expensive to evaluate and propose a solution to this problem using an Augmented Lagrangian (AL) with Gaussian Process (GP) emulation and an Expected Improvement criterion \citep{JonSchWel1998} for constrained minimization.  As pointed out by \cite{GraGraLeD2015}, statistical optimization is a massively expanding field and their paper addresses a significant problem within that field.  The majority of this discussion is focused on two particular choices in modeling, the first is the interesting idea of removing the maximum from the third term in  the AL and the second is related to addressing correlated constraint functions.  Section~\ref{sect:nomax} presents results from some numerical experiments and discusses the theoretical implications of dropping the maximum from the AL.  Section~\ref{sect:corrout} considers the potential impact on the authors' proposed expected improvement criterion of accounting for correlation among emulated constraints.

\section{The `NoMax' approach}
\label{sect:nomax}

\cite{GraGraLeD2015} motivate their methods using the simple ``toy problem'' of Section 4.3, specifically $\min f(x)$ subject to $c_1(x) \leq 0$ and $c_2(x) \leq 0$, where
	$$f(x) = x_1 + x_2 \,, \quad c_1(x)=\frac{3}{2}-x_1-2x_2-\frac{1}{2}\sin\left(2\pi(x_1^2-2x_2)\right) \,, \quad c_2(x)=x_1^2+x_2^2-\frac{3}{2} \,.$$
We shall further explore this problem by considering two new objective functions. First, we consider the linear objective
\begin{center}
{\bf Version 1:} $f(x)=x_1-x_2$
\end{center}
where the solution $-1$ at the location $(0,1)$ lies on the boundary of the parameter space as opposed to the boundary defined by the constraints.  In the second version we consider a nonlinear objective
\begin{center}
\noindent {\bf Version 2:} $f(x)=\frac{1}{2} (x_1-0.6)^2+(x_2-0.6)^2$
\end{center}
where the solution $0$ at the point $(0.6, 0.6)$ lies in the interior of the constraint set.

\begin{figure}[!ht]
\centering
\epsfig{file=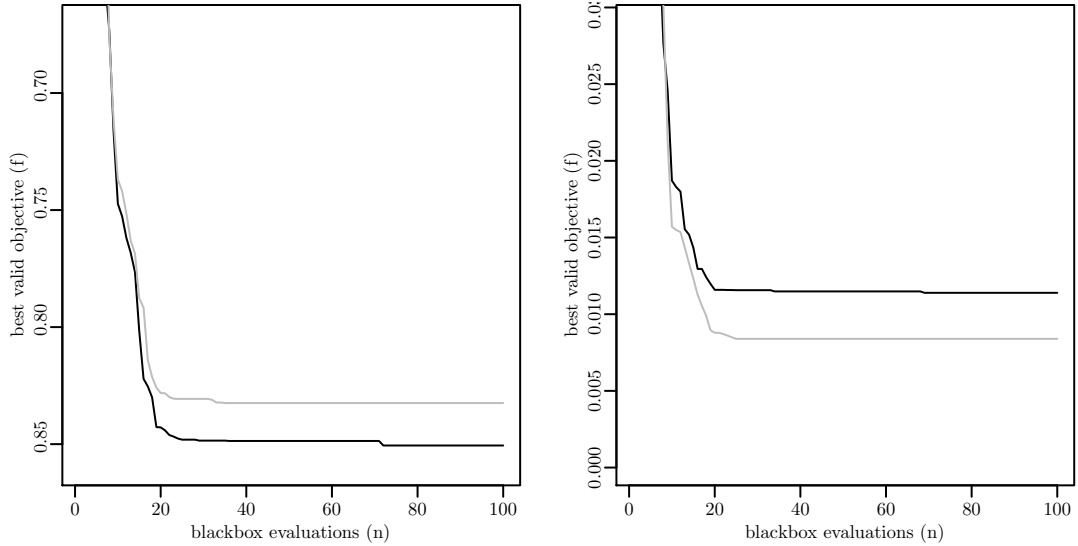}
\caption{Panel 1 (left) shows the results of solving Version 1 of the toy problem and Panel 2 (right) shows the results of solving Version 2 of the toy problem.  Black line is the solution of strategy {\bf NoMax} and the gray line is the solution of strategy {\bf WithMax}.}
\label{fig:nomax}
\end{figure}

For testing purposes we consider three different optimization strategies.  Strategy {\bf Random} is a naive approach of taking 100 uniform random samples and selecting the smallest function value that satisfies the constraints.  Strategy {\bf NoMax} uses the default setting of \verb!optim.auglag! in the \texttt{R} package \verb!laGP! \citep{Gra2014}, which removes the maximum from the AL.  Finally, strategy {\bf WithMax} uses \verb!optim.auglag!  with the option \verb!nomax=0!, which retains the maximum in the AL.  We follow the same setup as \cite{GraGraLeD2015} and consider 100 restarts of the algorithm.  The plot in Figure~\ref{fig:nomax} shows the results for both Versions 1 and 2 of the toy problem using strategies {\bf NoMax} and {\bf WithMax}.  Table~\ref{tab:res} shows the results of minimizing both versions of the toy problem using each of the three strategies.

\begin{table}[!ht]
\caption{Average value of the function after 100 function calls.}
\label{tab:res}
\begin{center}
\begin{tabular}{c|ccc}
 & {\bf Random} & {\bf NoMax} & {\bf WithMax}\\\hline
 {\bf $f(x) = x_1-x_2$} & -0.874 & -0.851 & -0.832 \\
 {\bf $f(x)=\frac{1}{2} (x_1-0.6)^2+(x_2-0.6)^2$} & 0.0014 &  0.0113&   0.0084
\end{tabular}
\end{center}
\end{table}

Figure \ref{fig:nomax} and Table \ref{tab:res} highlight two points.  First, at least in one case, the algorithm can be outperformed by a simple random search.  It is likely that this is partly reliant on the very low dimension of the toy problem, so should not be taken too seriously.  Second, when the test problem is nonlinear, the {\bf NoMax} is outperformed by the {\bf WithMax}.  This is actually expected, as explained below.

Recall from classic AL,
$$L_A(x; \lambda, \rho) := f(x) + \lambda^\top c(x) + \frac{1}{2 \rho}\sum_{i=1}^m (\max\{0, c_i(x)\})^2,$$
and the {\bf NoMax} version \citep{GraGraLeD2015},
$$L_A^{\tt nm}(x; \lambda, \rho) := f(x) + \lambda^\top c(x) + \frac{1}{2 \rho}\sum_{i=1}^m (c_i(x))^2.$$
The reason the AL is effective in practice is captured in the following Theorem.

\begin{thm} Let $f \in \mathcal{C}^2$ and $c_i \in \mathcal{C}^2$ for $i=1, 2, ..., m$.  Consider the problem of
    \begin{equation}\label{eq:min}
    \min \{ f(x) ~:~ c_i(x) \leq 0, ~i = 1, 2, ..., m, ~x \in \mathcal{B}\}.
    \end{equation}
Fix any $\rho>0$.  Suppose the Linear Independence Constraint Qualification (see \cite[Def 12.4]{NumOpt}) holds at $\tilde{x}$. Then $\tilde{x}$ is a critical point of Problem \eqref{eq:min} if and only if there are Lagrange multipliers $\tilde{\lambda} \geq 0$ such that $(\tilde{x}, \tilde{\lambda})$ is a critical point of the problem
    \begin{equation}\label{eq:augLA}
    \min_x \: \max_{\lambda \geq 0} L_A(x; \lambda, \rho).
    \end{equation}
\end{thm}

\proof ~It is immediately clear that Problem \eqref{eq:min} is exactly equivalent to
    \begin{equation}\label{eq:augF}
    \min \{ f(x)  + \frac{1}{2 \rho}\sum_{i=1}^m (\max\{0, c_i(x)\})^2 ~:~ c_i(x) \leq 0, ~i = 1, 2, ..., m, ~x \in \mathcal{B}\}.
    \end{equation}
Problem \eqref{eq:augLA} can now be viewed as the Lagrangian of Problem \eqref{eq:augF}.  The equivalence of Problem \eqref{eq:augF} and Problem \eqref{eq:augLA} now follows from standard Lagrange duality (see \cite[Thm. 12.1]{NumOpt} for example). {$\Box$} \medskip

\begin{figure}[!ht]
\centering
\includegraphics[width=1.25in]{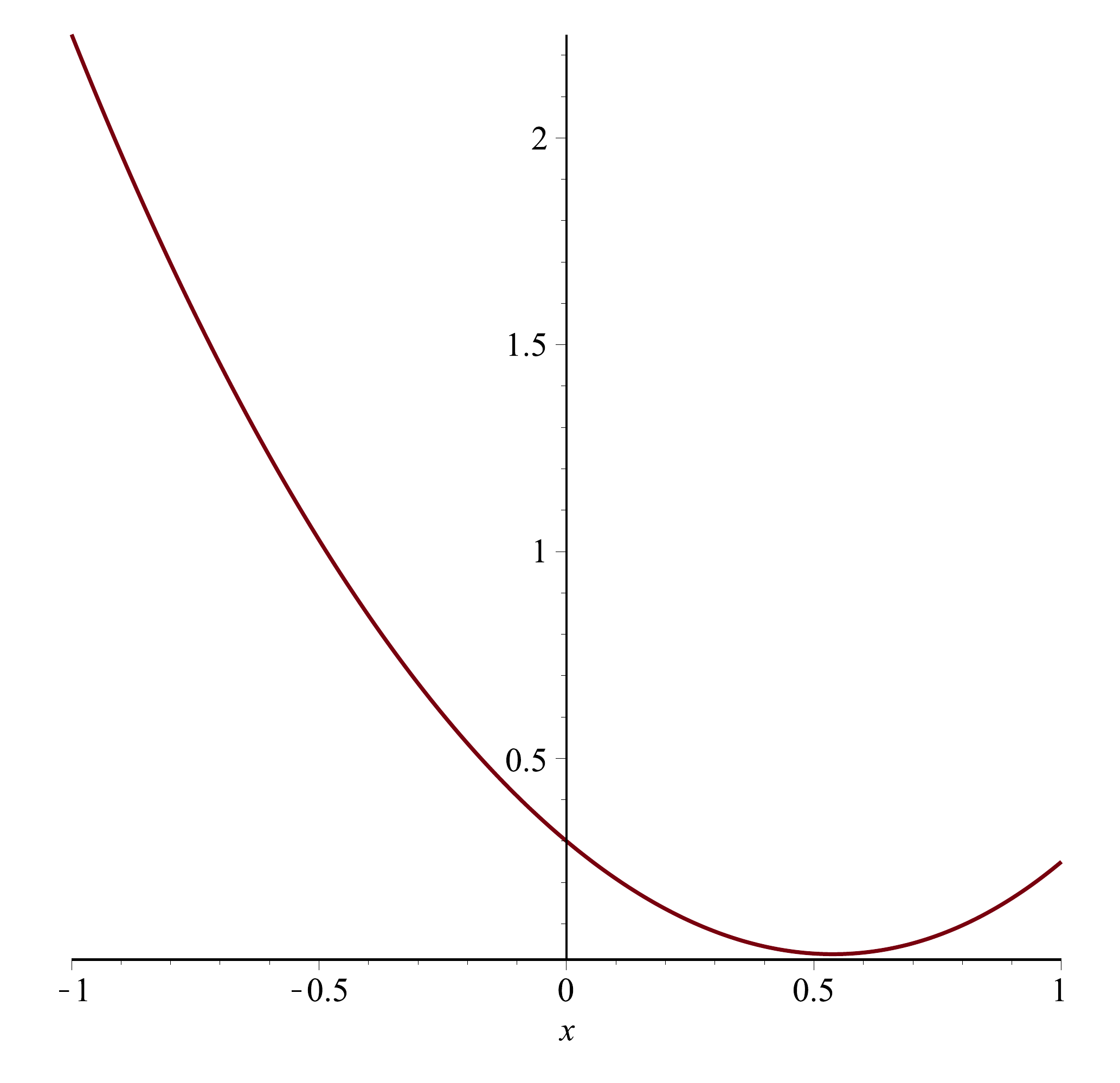} ~~\includegraphics[width=1.25in]{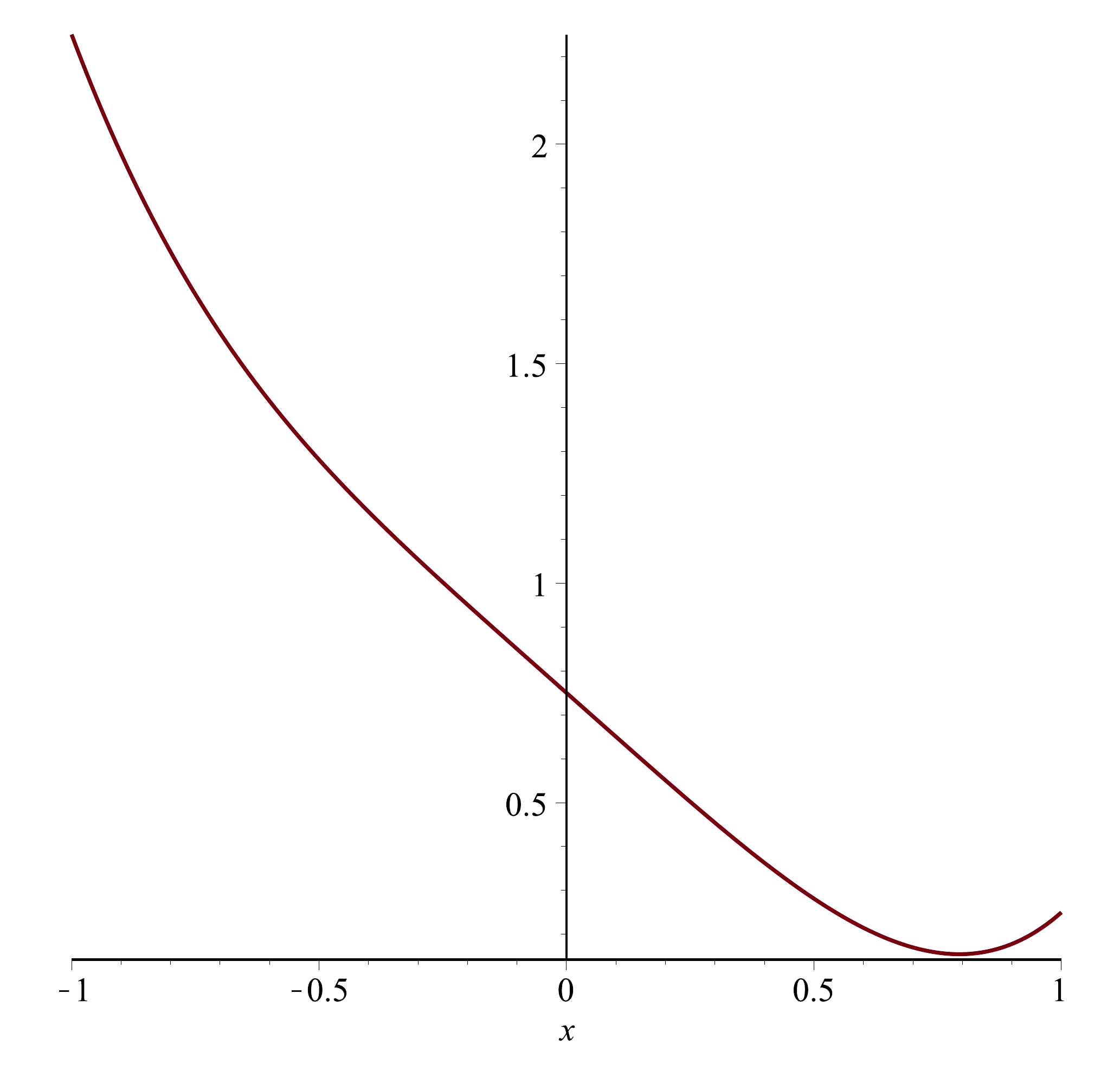}~~ \includegraphics[width=1.25in]{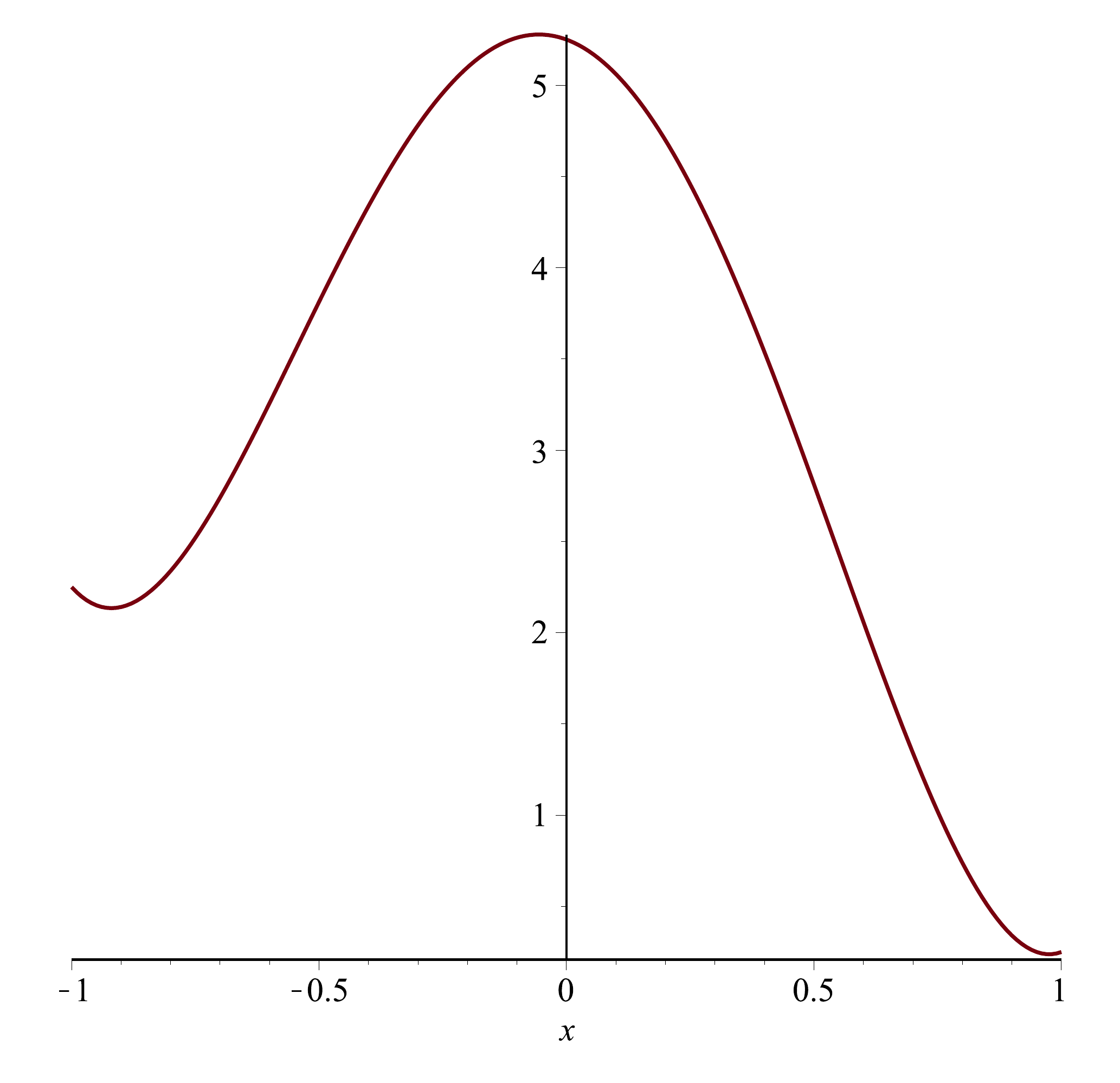}
\caption{Plots of $(x-0.5)^2 + \frac{1}{2 \rho}(x^2 - 1)^2$ for $\rho=10$ (left), $\rho=1$ (center) and $\rho=0.1$ (right).}
\label{fig:notzero}
\end{figure}

This theorem fails for the {\bf NoMax} version of the AL.  In order to illustrate we consider the one-dimensional problem
$$\min\{ (x-0.5)^2 : x^2 - 1 \leq 0 \}.$$
The minimizer of this problem is clearly unique and equal to $\tilde{x}=0.5$.  The {\bf NoMax} version of the AL for this problem is
$$\min_x \:\max_{\lambda \geq 0} L_A^{\tt nm}(x; \lambda, \rho)
    = \min_x \:\max_{\lambda \geq 0} \left\{ (x-0.5)^2 + \lambda (x^2-1) + \frac{1}{2 \rho}(x^2 - 1)^2 \right\}.$$
Examining $\max_{\lambda \geq 0} \{\lambda (x^2-1)\}$, we must have $(x^2-1) \leq 0$, or $\lambda$ will be driven to infinity and we are clearly not at the minimum with respect to $x$.  As such, we may reduce our AL subproblem to the one-dimensional problem
$$\min_{x} \left\{ (x-0.5)^2 + \frac{1}{2 \rho}(x^2 - 1)^2 ~:~ -1 \leq x \leq 1 \right\}.$$
In Figure \ref{fig:notzero} we provide several plots of $(x-0.5)^2 + \frac{1}{2 \rho}(x^2 - 1)^2$, using various values of $\rho$. The effect is clear.  As $\rho$ is decreased, the optimal point of the {\bf NoMax} AL is driven away from the true minimum and towards the boundary of the constraint set (note the algorithm drives $\rho$ to $0$ whenever stagnation is detected).  Moreover, regardless of how large $\rho$ becomes, $\min_x \:\max_{\lambda \geq 0} L_A^{\tt nm}(x; \lambda, \rho) $ will never return $\tilde{x}=0.5$ in this example.

Why then does the {\bf NoMax} version work for \cite{GraGraLeD2015}?  The reason lies in the use of the linear objective function.  In particular, if the objective function is linear (and nonzero), then the minimizer(s) will lie on the boundary of the constraint set.  Indeed, any point $\tilde{x}$ in the interior of the constraint set will have a nonzero gradient, so Fermat's Theorem confirms it cannot be a minimizer.

We point out these issues not as a flaw in the paper, but to highlight the importance of the assumption that the objective function be linear.
On the positive side, the assumption is not particularly demanding.  After all, if the objective function is nonlinear, one can always rewrite the optimization with a linear objective:
$$\min \{ f(x) : x \in \Omega\} = \min \{ r : r \geq f(x), x \in \Omega\}.$$

\begin{figure}[!ht]
\begin{center}
\epsfig{file=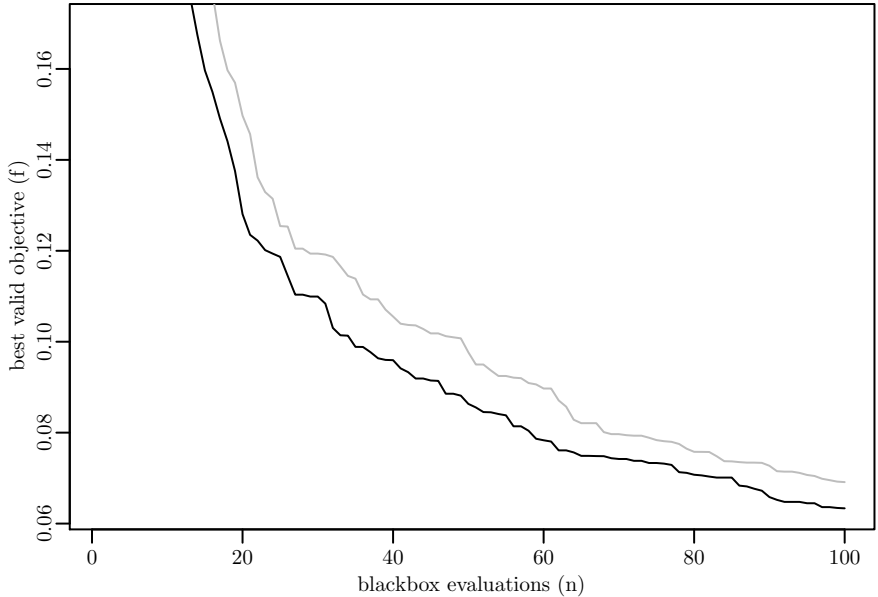}
\end{center}
\caption{Black line is the solution of {\bf NoMax} and the gray line is the solution of {\bf WithMax} for Version 3 of the toy problem.}
\label{fig:lincon}
\end{figure}

Following along these lines we consider Version 3 of the toy problem, a modification of Version 2 taking
$f(x)=x_3$ with constraints $c_1(x)$, $c_2(x)$ above and the additional constraint $c_3(x) \leq 0$, where
 $$c_3(x)=\frac{1}{2} (x_1-0.6)^2+(x_2-0.6)^2-x_3.$$ Again we use the same three fitting strategies as above. The plot in Figure~\ref{fig:lincon} shows the results using strategies {\bf NoMax} and {\bf WithMax}. The average function value using {\bf Random} was found to be 0.046 compared to average function values of 0.063 using {\bf NoMax} and 0.069 using {\bf WithMax}.

What is clear from each of these three simulations is that the constrained optimization algorithm proposed in \cite{GraGraLeD2015} shows promise but there is still work to be done in terms of testing and improving the method to deal with more complicated situations.  It should be noted that, if there are multiple constraint functions, with some active and some inactive at the minimizer, the {\bf NoMax} version may still present difficulties, as it may pull points towards inactive constraints.

\section{Correlated Outputs}
\label{sect:corrout}

We address the potential impact of explicitly modeling cross correlations among emulated outputs.  The authors considered only the independent output case, while acknowledging the potential for algorithmic efficiency improvements via such considerations.  As before, our brief investigation was carried out on the toy problem analyzed in Section 4.3 of \cite{GraGraLeD2015}.  The two constraint functions of this problem evaluated at a random sample of 1000 input pairs exhibited a correlation of approximately $-0.8$, suggesting the potential for improved emulation via modeling cross correlation in this application.  To keep algorithmic computation times at the same order as the independent output case, our correlated output model for constraints $c_1$ and $c_2$ was formulated as a simplified linear model of coregionalization (LMC),
\[
\begin{pmatrix}
c_1 (x) \\
c_2 (x)
\end{pmatrix} = A_n
\begin{pmatrix}
u_1 (x) \\
u_2 (x)
\end{pmatrix}
\]
where $u_1$ and $u_2$ are independent GPs and the fixed $2 \times 2$ coefficient matrix $A_n$ is calculated from the $n$ current observations of the two constraints placed in the $n \times 2$ matrix $C_n$ as follows:  Compute the ``economy size" singular value decomposition $C_n = U_n \Sigma_n V_n^T$ and set $A_n = V_n \Sigma_n$.  The two orthogonal columns of $U_n$ are transformed constraint observations modeled by the independent processes $u_1$ and $u_2$.  A more general LMC such as discussed in \cite{FriOakUrb2013} could be entertained as an improvement to our rudimentary approach, which will suffice for the current discussion.

\begin{figure}[!ht]
\centering
\includegraphics{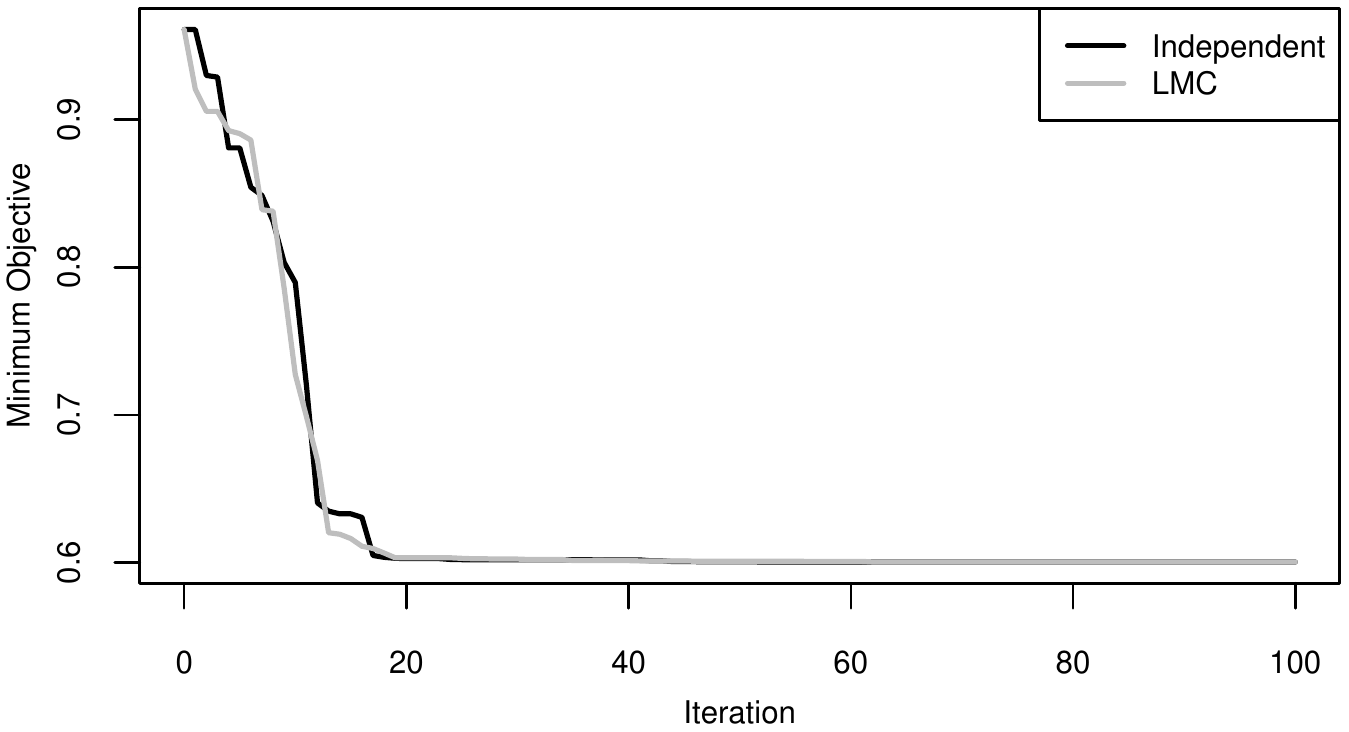}
\caption{Constrained minimum objective function value averaged over five 10-run initial designs for independent (black) and correlated (gray) output models.}
\label{fig:id1}
\end{figure}

Figure~\ref{fig:id1} shows the minimum objective value versus iteration for the scenario considered in \cite{GraGraLeD2015}, namely an initial design of 10 runs followed by 100 iterations of Algorithm 1.  These results are averaged over five initial Latin hypercube designs generated from the Matlab \texttt{lhsdesign} function with different seeds.  This small study suggests little to be gained from cross correlation modeling.

\begin{figure}[!ht]
\begin{center}
\includegraphics{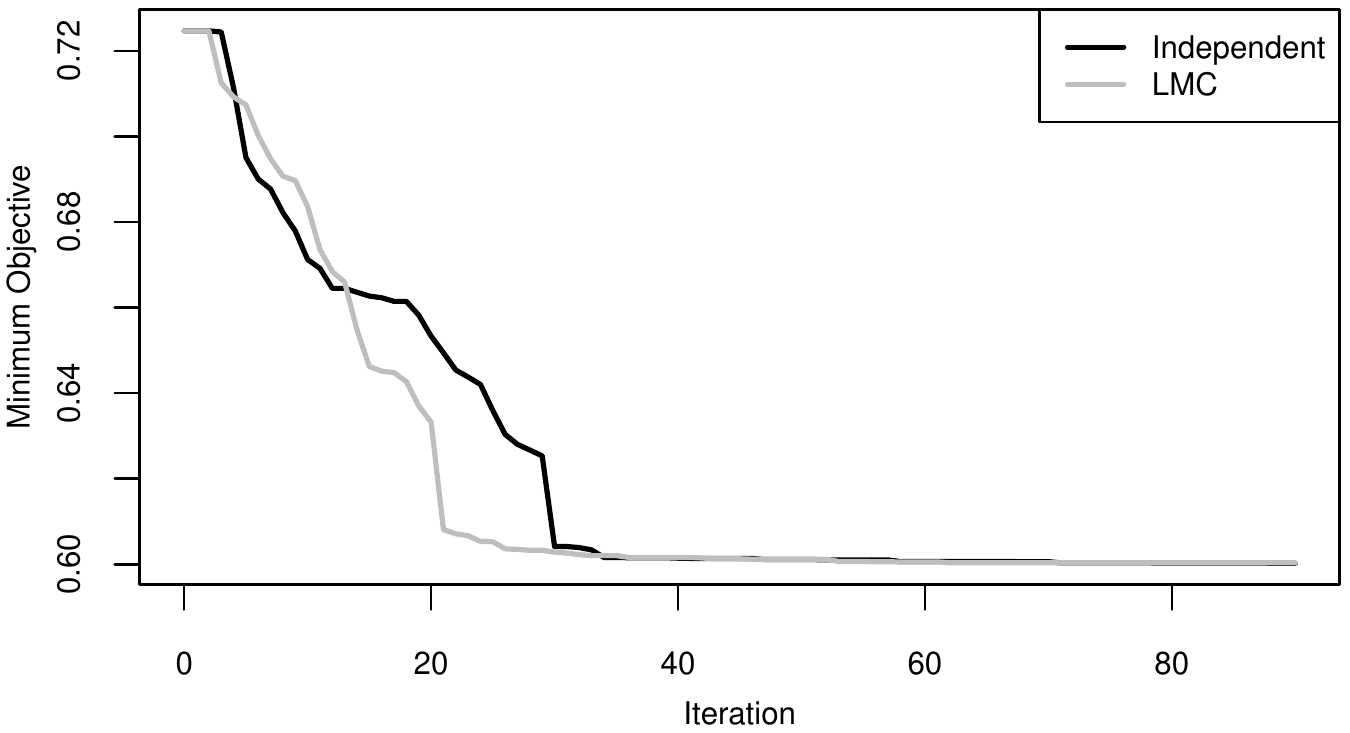}
\end{center}
\caption{Constrained minimum objective function value averaged over five 20-run initial designs for independent (black) and correlated (gray) output models.}
\label{fig:id2}
\end{figure}

It may be the case that the additional complication of modeling cross correlations is mitigated with larger initial designs at the expense of expected improvement iterations given a fixed budget.  Figure~\ref{fig:id2} shows the minimum objective value versus iteration averaged over five initial Latin hypercube designs of 20 runs generated as in the previous study followed by 90 iterations of Algorithm 1.  It appears in this case that cross correlation modeling could be beneficial for accelerating convergence to the constrained minimizer.  We emphasize again that the less sophisticated LMC model considered here requires roughly the same computational effort to implement as the independent output model.

Clearly no definitive conclusions can be reached regarding the efficacy of cross correlation modeling due to the small size and limited scope of the studies we conducted.  Nevertheless, such modeling appears to do no harm and is worthy of further study as a potentially useful alternative in expected improvement algorithms involving emulation of multiple correlated outputs.

\centerline{\bf ACKNOWLEDGEMENTS}
\renewcommand{\baselinestretch}{1.0}
The research of Loeppky was supported by Natural Sciences and Engineering Research Council of Canada Discovery Grant  ({RGPIN-2015-03895 }). The research of Hare was supported by Natural Sciences and Engineering Research Council of Canada Discovery Grant ({RGPIN-2013-355571}).  Los Alamos National Laboratory, an affirmative action/equal opportunity employer, is operated by Los Alamos National Security, LLC, for the National Nuclear Security Administration of the U.S. Department of Energy under contract DE-AC52-06NA25396.


\baselineskip 0.7\normalbaselineskip
\bibliographystyle{asa}
\bibliography{psn}

\end{document}